\documentclass{amsart} 
\usepackage{graphicx}
\usepackage{amsmath}
\usepackage{amssymb}
\usepackage{amsthm}
\newtheorem{thm}{Theorem}
\newtheorem{prop}{Proposition}

\newtheorem{lem}{Lemma} 
\newtheorem{definition}{Definition}
\newtheorem{rem}{Remark}

\begin{document}

\title[Remaining Case]{Radial Limits of Nonparametric PMC Surfaces with Intermediate Boundary Curvature}
\date{\today}
\author{Mozhgan ``Nora'' Entekhabi }
\address{Department of Mathematics \\ Florida A \& M University \\ Tallahassee, FL 32307}
\thanks{This research partially supported by NSF Award HRD-1824267.}
\email{mozhgan.entekhabi@famu.edu}
\author{Kirk E. Lancaster}
\address{Wichita, Kansas 67226}
\email{redwoodsrunner@gmail.com}

 
\def\Real{{\rm I\hspace{-0.2em}R}}
\def\Natural{{\rm I\hspace{-0.2em}N}}
\newcommand\myeq{\mathrel{\overset{\makebox[0pt]{\mbox{\normalfont\tiny\sffamily def}}}{=}}}

\begin{abstract}
We investigate the boundary behavior of the variational solution $f$  of a Dirichlet problem for a prescribed mean curvature equation 
in a domain $\Omega\subset\Real^{2}$  near a point $\mathcal{O}\in\partial\Omega$  under different assumptions about the curvature of 
$\partial\Omega$  on each side of $\mathcal{O}.$       
We prove that the radial limits at $\mathcal{O}$  of $f$  exist under different assumptions about the Dirichlet boundary data $\phi,$  
depending on the curvature properties of $\partial\Omega$  near $\mathcal{O}.$  
\end{abstract}

\maketitle

\section{Introduction}
Let $\Omega$  be a locally Lipschitz domain in $\Real^{2}$  and define $Nf = \nabla \cdot Tf = {\rm div}\left(Tf\right),$  where  $f\in C^{2}(\Omega)$  and 
$Tf= \frac{\nabla f}{\sqrt{1+\left|\nabla f\right|^{2}}}.$   
Let $H\in C^{1,\lambda}(\overline{\Omega})$  for some $\lambda\in (0,1)$  and satisfy the condition 
\[
\left|\int_{\Omega} H\eta \ dx\right| \le \frac{1}{2} \int_{\Omega} |D\eta|\ dx \ \ \ \ {\rm for \ all \ } \eta\in C^{1}_{0}(\Omega)
\]
(e.g.\cite[(16.60)]{GT}, \cite{Gui:78}). 
Here and throughout the paper, we adopt the sign convention that the curvature of $\Omega$  is nonnegative when $\Omega$  is convex.  
Consider the Dirichlet problem
\begin{eqnarray}
\label{eq:D}
Nf & = & 2H  \mbox{  \ in \ } \Omega    \\
f & = & \phi  \mbox{ \ on \ } \partial \Omega.
 \label{bc:D}
\end{eqnarray}  
Understanding the boundary behavior of a solution of (\ref{eq:D})-(\ref{bc:D}) has been the goal of many authors.

The geometry of $\Omega$  plays a critical role with regard to the existence of functions $f\in C^{2}(\Omega)\cap C^{0}(\overline{\Omega})$  which 
satisfy (\ref{eq:D}) and (\ref{bc:D})  (i.e. classical solutions of (\ref{eq:D})-(\ref{bc:D})). 
For some choices of domain $\Omega$  and boundary data $\phi,$  no classical solution of (\ref{eq:D})-(\ref{bc:D}) exists; when $H\equiv 0,$  much of the history 
(up to 1985) of this topic can be found in Nitsche's book \cite{NitscheBook} (e.g. \S 285, 403--418) and, for general $H,$  one might consult \cite{Ser:69}.  
(Appropriate ``smallness of $\phi$'' conditions can imply the existence of classical solutions when $\Omega$  is not convex in the $H\equiv 0$  case 
(e.g. \cite[\S 285 \& \S 412]{NitscheBook} and  \cite{Korn,RT,Williams1984,Williams1985}) or when $\partial\Omega$  does not satisfy appropriate curvature conditions 
in the general case   (e.g. \cite{Bergner,HayasidaNakatani,Lau}); however see \cite[\S 411]{NitscheBook}.) 
Different notions of ``generalized'' solutions of (\ref{eq:D})-(\ref{bc:D}) exist, such as Perron solutions (e.g. \cite{GT},\cite[\S 416]{NitscheBook})
and variational solutions (e.g. \cite{FinnBook},\cite[\S 417-418]{NitscheBook}); we shall focus on variational solutions.

The most extreme case (for locally Lipschitz domains in the plane) occurs when $\partial\Omega$  has a corner (or corners) and understanding   
the boundary behavior of solutions of (\ref{eq:D})-(\ref{bc:D}) near a corner is best investigated by understanding the radial limits of $f$  at the corner. 
The existence of radial limits when $H\equiv 0$  was established in \cite{Lan1985}  (see also \cite{EL1986A,EL1989,Lan1988}) and this was extended to general $H$ 
in \cite{EL1986B}  (see also \cite{LS1,NoraKirk1}).  

\begin{figure}[htb]
\centerline{
\includegraphics[width=3in]{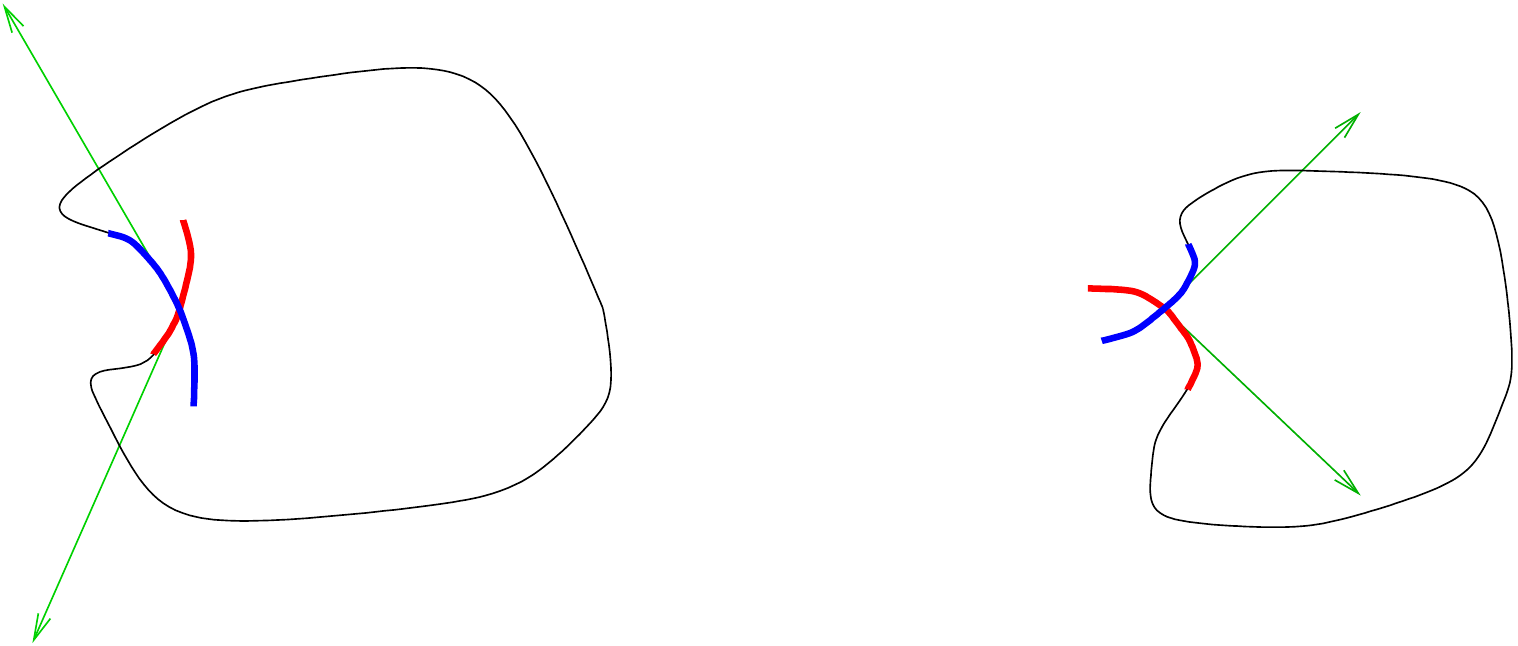}}
\caption{$\Sigma^{\pm}$  when $\beta-\alpha > \pi$  (left) \hspace{3mm}  $\Sigma^{\pm}$  when $\beta-\alpha < \pi$  (right)  \label{Fig0}}
\end{figure}

Let us assume that $\mathcal{O}=(0,0)\in\partial\Omega$  and there exist $\delta_{0}>0$  and $\alpha,\beta\in (-\pi,\pi)$  with $\alpha<\beta$  such that 
$B_{\delta_{0}}(\mathcal{O})\cap\partial\Omega\setminus\{\mathcal{O}\}$  consists of two components, $\partial^{-}\Omega$  and $\partial^{+}\Omega,$  
which are smooth (i.e. $C^{2,\lambda}$  for some $\lambda\in (0,1)$) curves, the rays $\theta=\alpha$  and $\theta=\beta$  are tangent rays to $\partial\Omega$  
at $\mathcal{O},$   $\partial\Omega$  has a corner at $\mathcal{O}$  of size $\beta-\alpha\in (0,2\pi)$  and 
\[
\{r(\cos\theta,\sin\theta) : 0<r<\epsilon(\theta), \alpha<\theta<\beta \} \subset \Omega\cap B_{\delta_{0}}(\mathcal{O}) 
\]
for some function $\epsilon(\cdot):(\alpha,\beta)\to (0,\delta_{0});$  here $(r,\theta)$  represents polar coordinates about $\mathcal{O}$  and   
$B_{\delta}(\mathcal{O})=\{{\bf x}\in\Real^{2} : |{\bf x}-\mathcal{O}|<\delta\}.$   
We assume $\partial^{-}\Omega$  is tangent to the ray $\theta=\alpha,$   $\partial^{+}\Omega$  is tangent to the ray $\theta=\beta$   at $\mathcal{O},$  
$\partial^{-}\Omega$  is an (open) subset of a $C^{2,\lambda}$-curve $\Sigma^{-}$  which contains $\mathcal{O}$  as an interior point and 
$\partial^{+}\Omega$  is an (open) subset of a $C^{2,\lambda}$-curve $\Sigma^{+}$  which contains $\mathcal{O}$  as an interior point; 
if $\beta-\alpha=\pi,$  we assume $\Sigma^{-}=\Sigma^{+}$   (see Figure \ref{Fig0}).

Let $f\in BV(\Omega)\cap C^{2}(\Omega)$  minimize the functional 
\begin{equation}
\label{EQ1}
J(h)=\int_{\Omega} \sqrt{1+|Dh|^{2}} + \int_{\Omega}  2H h d{\bf x}  + \int_{\partial\Omega} |u-\phi| dH_{1}   
\end{equation}
for $h\in BV(\Omega),$  so that $f$  is the variational solution of (\ref{eq:D})-(\ref{bc:D}).  
Let $Rf(\theta)$  denote the radial limit of $f$  at $\mathcal{O}$  in the direction $\theta\in (\alpha,\beta),$  
\[
Rf(\theta)=\lim_{r\downarrow 0} f(r\cos\theta,r\sin\theta),
\]
and set $Rf(\alpha)=\lim_{\partial^{-}\Omega\ni {\bf x}\to \mathcal{O}} f^{*}({\bf x})$  and 
$Rf(\beta)=\lim_{\partial^{+}\Omega\ni {\bf x}\to \mathcal{O}} f^{*}({\bf x})$  when these limits exist, where 
$f^{*}$  denotes the trace of $f$  on $\partial\Omega.$  
In \cite{NoraKirk1} (together with \cite{EchartLancaster1}), the following two results were proven.

\begin{prop} (\cite[Theorem 1]{NoraKirk1}; \cite{EchartLancaster1})
\label{CONCLUSION}
Let $f\in C^{2}(\Omega)\cap L^{\infty}(\Omega)$  satisfy (\ref{eq:D}) and suppose $\beta-\alpha>\pi.$
Then for each $\theta\in (\alpha,\beta),$  $Rf(\theta)$  exists and $Rf(\cdot)$  is a continuous function on $(\alpha, \beta)$ which behaves in one 
of the following ways:
\begin{itemize}
\item[(i)] $Rf$  is a constant function  and all nontangential limits of $f$  at $\mathcal{O}$  exist.  
\item[(ii)] There exist $\alpha_{1}, \alpha_{2}\in [\alpha,\beta]$  with $\alpha_{1}<\alpha_{2}$  such that 
\[
Rf(\theta) \ \ {\rm is} 
\left\{ \begin{array}{ccc} 
{\rm constant} & {\rm for} & \alpha<\theta\le\alpha_{1}\\ 
{\rm strictly \ monotonic} & {\rm for} & \alpha_{1}\le\theta\le\alpha_{2}\\
{\rm constant} & {\rm for} & \alpha_{2}\le\theta<\beta. \\ 
\end{array}
\right.
\]
\item[(iii)] There exist $\alpha_{1},$  $\alpha_{2}$  and $\theta_{0}$  with 
$\alpha \leq \alpha_{1} < \theta_{0} < \theta_{0}+\pi < \alpha_{2} \leq \beta$  such that  
\[
Rf(\theta) \ \ {\rm is} 
\left\{ \begin{array}{ccc} 
{\rm constant} & {\rm for} & \alpha<\theta\le\alpha_{1}\\ 
{\rm strictly \ increasing \ } & {\rm for} & \alpha_{1}\le\theta\le\theta_{0}\\
{\rm constant} & {\rm for} & \theta_{0}\le \theta\le\theta_{0}+\pi\\ 
{\rm strictly \ decreasing \ } & {\rm for} & \theta_{0}+\pi\le\theta\le\alpha_{2}\\
{\rm constant} & {\rm for} & \alpha_{2}\le\theta<\beta. \\ 
\end{array}
\right.
\]
\item[(iv)] There exist $\alpha_{1},$  $\alpha_{2}$  and $\theta_{0}$  with 
$\alpha \leq \alpha_{1} < \theta_{0} < \theta_{0}+\pi < \alpha_{2} \leq \beta$  such that  
\[
Rf(\theta) \ \ {\rm is} 
\left\{ \begin{array}{ccc} 
{\rm constant} & {\rm for} & \alpha<\theta\le\alpha_{1}\\ 
{\rm strictly \ decreasing \ } & {\rm for} & \alpha_{1}\le\theta\le\theta_{0}\\
{\rm constant} & {\rm for} & \theta_{0}\le \theta\le\theta_{0}+\pi\\ 
{\rm strictly \ increasing \ } & {\rm for} & \theta_{0}+\pi\le\theta\le\alpha_{2}\\
{\rm constant} & {\rm for} & \alpha_{2}\le\theta<\beta. \\ 
\end{array}
\right.
\]
\end{itemize}
 \end{prop}

\begin{prop}  (\cite[Theorem 2]{NoraKirk1}; \cite{EchartLancaster1})
\label{DEMONSTRATION}
Let $f\in C^{2}(\Omega)\cap L^{\infty}(\Omega)$  satisfy (\ref{eq:D}) and suppose $m=\lim_{\partial^{-}\Omega\ni {\bf x}\to \mathcal{O} } f\left({\bf x}\right)$  
exists.
Then for each $\theta\in (\alpha,\beta),$   $Rf(\theta)$  exists and $Rf(\cdot)$  is a continuous function 
on $[\alpha, \beta),$  where $Rf(\alpha) \myeq m.$  
If  $\beta-\alpha \le\pi,$  $Rf$  can behave as in (i) or (ii) in Proposition \ref{CONCLUSION}.  
If $\beta-\alpha >\pi,$  $Rf$  can behave as in (i), (ii), (iii) or (iv) in Proposition \ref{CONCLUSION}. 
\end{prop}

The necessity of assuming the existence of $\lim_{\partial^{-}\Omega\ni {\bf x}\to \mathcal{O} } f\left({\bf x}\right)$  
when $\beta-\alpha \le\pi$  in Proposition \ref{DEMONSTRATION} is illustrated by the use of the ``gliding hump'' construction in \cite{Lan:89} 
and \cite[Theorem 3]{LS1}, where examples of $\Omega$  (with $\beta-\alpha=\pi$), $\phi$  and $H$  are presented such that $f$  is discontinuous at ${\mathcal O}$  
and none of the radial limits of $f$  at ${\mathcal O}$  exist.  
This same construction can be used to obtain examples of $\Omega,$  $\phi$  and $H$  with $0<\beta-\alpha<\pi$  such that $f$  is discontinuous at ${\mathcal O}$  
and none of the radial limits of $f$  at ${\mathcal O}$  exist.  

Thus the size $\beta-\alpha$  of the angle made by $\partial\Omega$  at the corner $\mathcal{O}$  is a primary determinant of the existence of radial limits 
$Rf(\theta),$  as Proposition \ref{CONCLUSION} shows that $Rf(\theta)$  exists for $\alpha<\theta<\beta$  for any solution $f\in C^{2}(\Omega)\cap L^{\infty}(\Omega)$  
of (\ref{eq:D}) when $\beta-\alpha>\pi,$  without regard to the behavior of (the trace of)  $f$  on $\partial\Omega.$  
An important question is ``Does a  solution of (\ref{eq:D})-(\ref{bc:D}) actually satisfy the boundary condition (\ref{bc:D})   
near a specific point (e.g. a corner)  $\mathcal{O}\in\partial\Omega$?''  
In particular, \cite{Lan1985}  and \cite{EL1986B} require the answer to this question to be ``yes.'' 
The answer to this question depends largely on the curvature of $\partial\Omega$  on each side of $\mathcal{O}$  and this is a secondary determinant of 
the existence of radial limits $Rf(\theta)$  when $\beta-\alpha \le\pi.$ 

To illustrate the importance of curvature conditions on the possible behaviors of solutions of (\ref{eq:D})-(\ref{bc:D}), 
suppose $\delta>0,$   $B_{\delta}(\mathcal{O})\cap\partial\Omega$  is smooth (so $\beta-\alpha=\pi$), $\kappa({\bf x})$  is the curvature of $\partial\Omega$  
at ${\bf x}\in B_{\delta}(\mathcal{O})\cap\partial\Omega,$  $\phi\in L^{\infty}(\partial\Omega)$  and $f$  is the variational solution of (\ref{eq:D})-(\ref{bc:D}). 
In \cite[Theorem 1.1]{NoraKirk3}, the authors proved the existence of radial limits $Rf(\theta),$  $\theta\in [\alpha,\beta],$  when 
$\kappa({\bf x})<-2|H({\bf x})|$  for ${\bf x}\in B_{\delta}(\mathcal{O})\cap\partial\Omega,$   without regard to the behavior of $\phi$  on $\partial\Omega.$ 

Our goals here are, in Theorem \ref{Thm1}, to extend the results in \cite{NoraKirk3} to the ``remaining case'' noted there in which $\partial\Omega$  is smooth and $\kappa$  satisfies 
\begin{equation}
\label{WWW}
-2|H({\bf x})|\le \kappa({\bf x})<2|H({\bf x})| \ \ \ {\rm for} \  {\bf x}\in B_{\delta}(\mathcal{O})\cap\partial\Omega 
\end{equation}
and, in Theorem \ref{Thm2}, to extend the results in \cite{NoraKirk3} to actual corners (i.e. $\beta-\alpha\neq\pi$).

\section{Theorems}

\begin{thm}
\label{Thm1}
Let $f$  be the variational solution of (\ref{eq:D})-(\ref{bc:D}).
Suppose $\Gamma\subset \partial\Omega$  is a $C^{2,\lambda}$  (open) curve  for some $\lambda\in (0,1),$  $\mathcal{O}\in \Gamma,$    
$H$  is non-negative or non-positive in a neighborhood of $\mathcal{O}$  and $\kappa(\mathcal{O})<2|H(\mathcal{\mathcal{O}})|.$  
Then $Rf(\theta)$  exists for each $\theta\in (\alpha,\beta),$   $Rf\in C^{0}((\alpha,\beta))$  and $Rf$  can behave as in (i) or (ii) in Proposition \ref{CONCLUSION}.
Further, if $\kappa(\mathcal{O})<-2|H(\mathcal{\mathcal{O}})|,$  then $Rf(\alpha)$  and $Rf(\beta)$  both exist, $Rf\in C^{0}([\alpha,\beta]),$  and, 
in case (i) in Proposition \ref{CONCLUSION}, $f\in C^{0}(\Omega\cup\{\mathcal{O}\}).$  
\end{thm}
\vspace{3mm}

\begin{thm}
\label{Thm2}
Suppose $H$  is non-negative or non-positive in a neighborhood of $\mathcal{O},$  
\begin{equation}
\label{curve2}
\limsup_{\partial^{\pm}\Omega\ni {\bf x}\to \mathcal{O}}\left(\kappa({\bf x})-2|H({\bf x})|\right)<0
\end{equation}
and, if  $\beta-\alpha < \pi,$   $m=\lim_{\partial^{-}\Omega\ni {\bf x}\to \mathcal{O} } f\left({\bf x}\right)$  exists.      
Then $Rf(\theta)$  exists for each $\theta\in (\alpha,\beta)$  and  $Rf\in C^{0}(\alpha,\beta).$         
\begin{itemize}
\item[(a)] If  $\beta-\alpha \le \pi,$  $Rf$  behaves as in (i) or (ii) in Proposition \ref{CONCLUSION}.  
\item[(b)] If $\beta-\alpha >\pi,$  $Rf$  behaves as in (i), (ii), (iii) or (iv) in Proposition \ref{CONCLUSION}. 
\end{itemize}

\noindent If, in addition,  $\limsup_{\partial^{\pm}\Omega\ni {\bf x}\to \mathcal{O}}\left(\kappa({\bf x})+2|H({\bf x})|\right)<0,$   
then $Rf(\alpha)$  and $Rf(\beta)$  both exist, $Rf\in C^{0}([\alpha,\beta]),$  and, in case (i) in Proposition \ref{CONCLUSION}, 
$f\in C^{0}(\Omega\cup\{\mathcal{O}\}).$  
\end{thm}
\vspace{3mm}

\noindent We note that the ``gliding hump'' construction (which depends on the existence of classical solutions of (\ref{eq:D})-(\ref{bc:D})) 
cannot be successfully used when $\beta-\alpha>\pi;$   however it remains an open question if radial limits of $f$  always exist without regard to the behavior 
of (the trace of)  $f$  on $\partial\Omega$ when (\ref{curve2}) holds and $\beta-\alpha<\pi$  (see Remark \ref{Why}).

\section{Proofs}

Let  $Q$  be the operator on $C^{2}(\Omega)$  given by 
\begin{equation}
\label{Q}
Qf({\bf x}) \myeq Nf({\bf x}) - 2H({\bf x}), \ \ \ \ {\bf x}\in\Omega. 
\end{equation}
Let $\nu$  be the exterior unit normal to $\partial\Omega,$  defined almost everywhere on $\partial\Omega.$    
At every point ${\bf y}\in\partial\Omega$  for which $\partial\Omega$  is a $C^{1}$  curve in a neighborhood of ${\bf y},$  
$\hat\nu$   denotes a continuous extension  of $\nu$  to a neighborhood of ${\bf y}.$ 
Finally we adopt the convention used in \cite[p. 178]{CF:74a} with regard to the meaning of phrases like 
``$T\psi({\bf y})\cdot\nu({\bf y})=1$  at a point ${\bf y}\in \partial\Omega$'' and the notation, definitions and conventions used in \cite{NoraKirk3}, 
including upper and lower Bernstein pairs $\left(U^{\pm},\psi^{\pm}\right),$  which we quote below.
\vspace{3mm}

\begin{definition}
\label{Bernstein+} 
Given a locally Lipschitz domain $\Omega,$ an \underline{\bf upper Bernstein pair} $\left(U^{+},\psi^{+}\right)$  for a curve $\Gamma\subset \partial\Omega$  
and a function $H$  in (\ref{Q}) is a domain $U^{+}$  and a function $\psi^{+}\in C^{2}(U^{+})\cap C^{0}\left(\overline{U^{+}}\right)$  such that 
$\Gamma\subset\partial U^{+},$  $\nu$  is the exterior unit normal to $\partial U^{+}$  at each point of $\Gamma$  
(i.e. $U^{+}$  and $\Omega$  lie on the same side of $\Gamma$), $Q\psi^{+}\le 0$  in $U^{+},$  and $T\psi^{+}\cdot\nu=1$  almost everywhere on $\Gamma$  in the same sense as in  \cite{CF:74a}; that is, for almost every ${\bf y}\in\Gamma,$ 
\begin{equation}
\label{zero1}
\lim_{U^{+}\ni {\bf x}\to {\bf y}} \frac{\nabla \psi^{+}({\bf x})\cdot \hat\nu({\bf x})}{\sqrt{1+|\nabla \psi^{+}({\bf x})|^{2}}} = 1.
\end{equation} 
\end{definition}

\begin{definition}
\label{Bernstein-} 
Given a domain $\Omega$  as above, a \underline {\bf lower Bernstein pair} $\left(U^{-},\psi^{-}\right)$  for a curve $\Gamma\subset \partial\Omega$  
and a function $H$  in (\ref{Q}) is a domain $U^{-}$  and a function $\psi^{-}\in C^{2}(U^{-})\cap C^{0}\left(\overline{U^{-}}\right)$  
such that $\Gamma\subset\partial U^{-},$  $\nu$  is the exterior unit normal to $\partial U^{-}$  at each point of $\Gamma$  
(i.e. $U^{-}$  and $\Omega$  lie on the same side of $\Gamma$), 
$Q\psi^{-}\ge 0$  in $U^{-},$  and $T\psi^{-}\cdot\nu=-1$  almost everywhere on $\Gamma$  (in the same sense as above).
\end{definition}

\noindent The argument which establishes \cite[Corollary 14.13]{GT}, together with boundary regularity results (e.g. \cite{Bour,Lin}), are noted in 
\cite[Remark 1]{NoraKirk3} and imply the following 

\begin{lem}
\label{lem1}
Suppose $\Delta$  is a $C^{2,\lambda}$  domain in $\Real^{2}$  for some $\lambda\in (0,1),$   ${\bf y}\in\partial\Omega$  and 
$\Lambda({\bf y})< 2|H({\bf y})|,$  where $\Lambda({\bf y})$  denotes the curvature of $\partial\Delta$  at  ${\bf y}.$  
If $H$  is non-negative  in $U\cap\Omega$  for some neighborhood $U$  of ${\bf y},$  then there exist $\tau>0$  and an upper Bernstein pair 
$\left(U^{+},\psi^{+}\right)$  for $(\Gamma,H),$  where $\Gamma= B_{\tau}({\bf y})\cap\partial\Omega$  and $U^{+}=B_{\tau}({\bf y})\cap \Omega.$    
If $H$  is non-positive  in $U\cap\Omega$  for some neighborhood $U$  of ${\bf y},$  then there exist $\tau>0$  and a lower Bernstein pair 
$\left(U^{-},\psi^{-}\right)$  for $(\Gamma,H),$  where $\Gamma= B_{\tau}({\bf y})\cap\partial\Omega$  and $U^{-}=B_{\tau}({\bf y})\cap \Omega.$  
\end{lem}
\vspace{3mm} 

\noindent {\bf Proof of Theorem \ref{Thm1}:}  
We note, as in \cite{LS1}, that the conclusion of Theorem \ref{Thm1} is a local one and we may assume $\Omega$  is a bounded domain. 
The claims in the last sentence of the theorem follow from \cite[Theorem 1.1]{NoraKirk3}.  
We may assume that $f\in C^{0}(\overline{\Omega}\setminus \{\mathcal{O}\})$  (i.e. $f\in C^{2}(\Omega)$  and, if necessary, we could replace $\Omega$  
by a set $U\subset\Omega$  such that $\partial U \cap \partial\Omega = \{\mathcal{O}\},$  $\partial U$  has the same tangent rays at $\mathcal{O}$  
as does $\partial\Omega$  and the curvature $\kappa^{*}$  of $\partial U$  satisfies $\kappa^{*}(\mathcal{O})<2|H(\mathcal{\mathcal{O}})|$).  

Let $z_{1}= \liminf_{\Omega\ni {\bf x}\to \mathcal{O}} f({\bf x})$  and $z_{2}= \limsup_{\Omega\ni {\bf x}\to \mathcal{O}} f({\bf x});$  
if $z_{1}=z_{2},$  then (i) of Proposition \ref{CONCLUSION} holds and thus we assume $z_{1}<z_{2}.$  
Set  $S_{0} = \{ ({\bf x},f({\bf x})) : {\bf x} \in \Omega \}.$
Since $f$  minimizes $J$  in (\ref{EQ1}),  we see that the area of $S_{0}$  is finite; let $M_{0}$  denote this area. 
For $\delta\in (0,1),$  set 
\[
p(\delta) = \sqrt{\frac{8\pi M_{0}}{\ln\left(\frac{1}{\delta}\right)}}.
\]
Let $E= \{ (u,v) : u^{2}+v^{2}<1 \}.$ 
As in \cite{EL1986A,LS1}, there is a parametric description of the surface $S_{0},$  
\begin{equation}
\label{PARAMETRIC1}
Y(u,v) = (a(u,v),b(u,v),c(u,v)) \in C^{2}(E:{\Real}^{3}), 
\end{equation}
which has the following properties:
\begin{itemize}
\item[$\left(a_{1}\right)$]  $Y$ is a diffeomorphism of $E$ onto $S_{0}$.
\item[$\left(a_{2}\right)$]  Set $G(u,v)=(a(u,v),b(u,v)),$  $(u,v)\in E.$   Then $G \in C^{0}(\overline{E} : {\Real}^{2}).$  
\item[$\left(a_{3}\right)$]  Set $\sigma(\mathcal{O})=G^{-1}\left(\partial \Omega\setminus \{\mathcal{O}\}\right);$  
then $\sigma(\mathcal{O})$ is a connected (open) arc of $\partial E$  and $Y$ maps $\sigma(\mathcal{O})$  onto $\partial \Omega\setminus \{\mathcal{O}\}.$  
We may assume the endpoints of $\sigma(\mathcal{O})$  are ${\bf o}_{1}$  and ${\bf o}_{2}.$  
(Note that ${\bf o}_{1}$  and ${\bf o}_{2}$  are not assumed to be distinct.)
\item[$\left(a_{4}\right)$]  $Y$ is conformal on $E$: $Y_{u} \cdot Y_{v} = 0, Y_{u}\cdot Y_{u} = Y_{v}\cdot Y_{v}$   on $E$.
\item[$\left(a_{5}\right)$]  $\triangle Y := Y_{uu} + Y_{vv} = 2H\left(Y\right)Y_{u} \times Y_{v}$  on $E$.
\end{itemize}
Let $\zeta(\mathcal{O})=\partial E\setminus\sigma(\mathcal{O});$  then $G(\zeta(\mathcal{O}))=\{\mathcal{O}\}$  and 
${\bf o}_{1}$  and ${\bf o}_{2}$  are the endpoints of $\zeta(\mathcal{O}).$   
\vspace{2mm}

Suppose first that ${\bf o}_{1}\neq {\bf o}_{2}.$  
From the Courant-Lebesgue Lemma (e.g. Lemma $3.1$ in \cite{Cour:50}), we see that there exists 
$\rho=\rho(\delta,{\bf w})\in \left(\delta,\sqrt{\delta}\right)$  such that the arclength $l_{\rho}=l_{\rho(\delta,{\bf w})}$  of 
$Y(C_{\rho(\delta,{\bf w})}({\bf w}))$  is less than $p(\delta),$  for each $\delta\in (0,1)$  and ${\bf w}\in \partial E;$ 
here $C_{r}({\bf w})=\{(u,v)\in E : |(u,v)-{\bf w}|=r \}.$  
Set $E_{r}({\bf w})=\{(u,v)\in E : |(u,v)-{\bf w}|<r\},$  $E'_{r}({\bf w})=G(E_{r}({\bf w}))$  and  $C'_{r}({\bf w})=G(C_{r}({\bf w})).$  
Choose $\delta_{1}>0$  such that $2\sqrt{\delta_{1}}<|{\bf o}_{1}-{\bf o}_{2}|.$  
Let ${\bf w}_{0}\in \zeta(\mathcal{O})$  be the ``midpoint'' of ${\bf o}_{1}$  and ${\bf o}_{2},$  so that 
$\sqrt{\delta_{1}}<|{\bf w}_{0}-{\bf o}_{1}|=|{\bf w}_{0}-{\bf o}_{2}|.$  
Set $\mathcal{C}=C_{\rho(\delta_{1},{\bf w}_{0})}'({\bf w}_{0});$
then  $\{({\bf x},f({\bf x})) : {\bf x}\in \mathcal{C} \}$  ($=Y(C_{\rho(\delta_{1},{\bf w}_{0})}({\bf w}_{0}))$)
is a curve of finite length $l_{\rho(\delta_{1},{\bf w}_{0})}$  with endpoints $(\mathcal{O},z_{a})$  and $(\mathcal{O},z_{b})$  for some 
$z_{a},z_{b}\in\Real.$  
Notice, in particular, that the graph of $f$  over $\mathcal{C}$  is either continuous at $\mathcal{O}$  (if $z_{a}=z_{b}$)  or has a jump 
discontinuity at $\mathcal{O}$ (if $z_{a}\neq z_{b}$).   

We may now argue as in \cite{Lan1988}. 
Let $\Omega_{0}=G(E_{\rho(\delta_{1},{\bf w}_{0})}({\bf w}_{0}))=E'_{\rho(\delta_{1},{\bf w}_{0})}({\bf w}_{0}),$  
so that $\partial\Omega_{0}=\mathcal{C}\cup \{\mathcal{O}\}.$   
From the Courant-Lebesgue Lemma and the general comparison principle (\cite[Theorem 5.1]{FinnBook}), 
we see that $Y$  is uniformly continuous on $E_{\rho(\delta_{1},{\bf w}_{0})}({\bf w}_{0})$  
and so extends to a continuous function on the closure of $E_{\rho(\delta_{1},{\bf w}_{0})}({\bf w}_{0}).$
From Steps 2, 4 and 5 of \cite{LS1} and with \cite{EchartLancaster1} replacing Step 3 of \cite{LS1}, we see that 
there exist $\alpha_{0},\beta_{0}\in [\alpha,\beta]$  with $\alpha_{0}<\beta_{0}$  such that 
\[
\{r(\cos\theta,\sin\theta) : 0<r<\epsilon_{0}(\theta), \alpha_{0}<\theta<\beta_{0} \} \subset \Omega_{0}\cap B_{\delta_{0}}(\mathcal{O}) 
\]
for some function $\epsilon_{0}(\cdot):(\alpha,\beta)\to (0,\delta_{0})$  and the radial limits $Rf(\theta)$  of $f$  at $\mathcal{O}$  exist 
for $\alpha_{0}\le\theta\le\beta_{0}.$  
Since $\partial\Omega$  is ($C^{2,\lambda}$) smooth near $\mathcal{O},$  we have $\beta-\alpha=\pi$  and so $\beta_{0}-\alpha_{0}\le\pi.$    
(We note that $z_{a}=z_{b}$  when ${\bf o}_{1}\neq {\bf o}_{2}$  and $\beta_{0}-\alpha_{0}\le\pi$  implies 
$f\in C^{0}(\overline{\Omega}),$  a contradiction, and so $z_{a}\neq z_{b}.$)    
The existence of $Rf(\cdot)$  on $(\alpha,\beta)$  now follows from two applications of \cite[Theorem 2]{NoraKirk1}, one in the domain 
$(r\cos\theta,r\sin\theta)\in\Omega : r>0, (\alpha_{0}+\beta_{0})/2<\theta<\beta\}$  and one in the domain 
$(r\cos\theta,r\sin\theta)\in\Omega : r>0, \alpha<\theta<(\alpha_{0}+\beta_{0})/2\}.$

Suppose second that ${\bf o}={\bf o}_{1}= {\bf o}_{2}$  and $\zeta(\mathcal{O})=\{{\bf o}\}.$   
Let us assume that $H$  is non-negative in a neighborhood of $\mathcal{O};$  here $H(Y(u,v))$  means $H(a(u,v),b(u,v)).$    
From Lemma \ref{lem1}, we see that an upper  Bernstein pair $\left(U^{+},\psi^{+}\right)$   for $(\Gamma_{1},H)$    exists, 
where $U^{+}=\Omega\cap B_{\tau}(\mathcal{O})$  and $\Gamma_{1}=\Gamma\cap B_{\tau}(\mathcal{O})$  for some $\tau>0;$    
let $q$  denote a modulus of continuity for  $\psi^{+}.$  
Then $T\psi^{+}\cdot \nu=+ 1$  (in the sense of \cite{CF:74a})  on $\Gamma_{1}$  and, for each $C\in\Real,$   
$Q(\psi^{+}+C)=Q(\psi^{+})\le 0$  on $\Omega\cap U^{+}$  or equivalently 
\begin{equation}
\label{Barrier}
N(\psi^{+}+C)({\bf x}) \le 2H({\bf x})=Nf({\bf x}) \ \ {\rm for} \ \ {\bf x}\in\Omega\cap U^{+}. 
\end{equation}

From the Courant-Lebesgue Lemma, we see that there exists 
$\rho=\rho(\delta,{\bf w})\in \left(\delta,\sqrt{\delta}\right)$  such that the arclength $l_{\rho}=l_{\rho(\delta,{\bf w})}$  of 
$Y(C_{\rho(\delta,{\bf w})}({\bf w}))$  is less than $p(\delta),$  for each $\delta\in (0,1)$  and ${\bf w}\in \partial E.$

Let us assume that $\delta\in (0,1)$  is small enough that $p(\delta)<\tau,$  so that $G({\bf w})\in U^{+}$  for each ${\bf w}\in E$  
with $|{\bf w}-{\bf o}|\le \sqrt{\delta}$  and $G({\bf w})\in \Gamma_{1}$  for each ${\bf w}\in \partial E$  with $|{\bf w}-{\bf o}|\le \sqrt{\delta}.$ 
Now $\psi^{+}-\psi^{+}({\bf x})\le q(p(\delta))$  in $E'_{\rho(\delta,{\bf o})}({\bf o})$  for any ${\bf x}\in E'_{\rho(\delta,{\bf o})}({\bf o})$  and   
the general comparison principle (\cite[Theorem 5.1]{FinnBook}) implies that if $U\subset E'_{\rho(\delta,{\bf o})}({\bf o})$  is an open set, then 
\begin{equation}
\label{wheat}
f \le \sup_{\Omega\cap\partial U}f+\psi^{+}-\inf_{\Omega\cap\partial U}\psi^{+}\le \sup_{\Omega\cap\partial U}f+q(p(\delta)) 
\ \ \ \ {\rm in} \ U.
\end{equation}
Set
\[
k(\delta)= \inf_{{\bf u}\in C_{\rho(\delta,{\bf o})}({\bf o})}c({\bf u}) = \inf_{ {\bf x}\in C'_{\rho(\delta,{\bf o})}({\bf o}) } f({\bf x}).
\]
Now $f\le k(\delta)+p(\delta)$  on $C'_{\rho(\delta,{\bf o})}({\bf o})$  and  $\psi^{+}-\inf_{C'_{\rho(\delta,{\bf o})}({\bf o})}\psi^{+}\le q(p(\delta))$  
in $E'_{\rho(\delta,{\bf o})}({\bf o})$     and so (\ref{wheat})  implies 
\[
f \le k(\delta)+p(\delta)+\psi^{+}-\inf_{C'_{\rho(\delta,{\bf o})}({\bf o})}\psi^{+}\le k(\delta)+p(\delta)+q(p(\delta)) 
\]
or
\begin{equation}
\label{rice}
\sup_{E'_{\rho(\delta,{\bf o})}({\bf o})} f  \le \inf_{C'_{\rho(\delta,{\bf o})}({\bf o})} f+p(\delta)+q(p(\delta)).
\end{equation}
Since $\sup_{E'_{\rho(\delta,{\bf o})}({\bf o}))} f\ge z_{2},$  
\begin{equation}
\label{oats}
\inf_{C'_{\rho(\delta,{\bf o})}({\bf o})} f \ge z_{2}-p(\delta)-q(p(\delta)) = z_{2}-o(\delta) \ \ \ 
{\rm for \ each} \ \ \delta>0.
\end{equation}
Let $z(\delta)=z_{2}-2p(\delta)-q(p(\delta))$  and 
\[
M(\delta)=\{{\bf x}\in \Omega\cap \overline{E'_{\rho(\delta,{\bf o})}({\bf o})} : f({\bf x})> z(\delta) \}.
\]
(Recall $f\in C^{0}(\overline{\Omega}\setminus \{\mathcal{O}\})$  and $c\in C^{0}(\overline{E}\setminus \{{\bf o}\}).$)   
Then for each $\delta\in (0,p^{-1}(\tau)),$  (\ref{oats}) implies $f\ge z_{2}-p(\delta)-q(p(\delta))>z(\delta)$  on $C'_{\rho(\delta,{\bf o})}({\bf o})$  and so 
\[
C'_{\rho(\delta,{\bf o})}({\bf o})\subset M(\delta) \ \ \ {\rm and} \ \ \ \mathcal{O}\in \overline{M(\delta)},
\] 
Let $V(\delta)$  denote the component of $M(\delta)$  which contains $C'_{\rho(\delta,{\bf o})}({\bf o}).$  
We claim that $\mathcal{O}\in \overline{V(\delta)}.$  
Suppose otherwise; then there is a curve $\mathcal{I}$  in $E'_{\rho(\delta,{\bf o})}({\bf o})$  (with endpoints ${\bf x}^{-}\in\partial^{-}\Omega$  and 
${\bf x}^{+}\in\partial^{+}\Omega$)  such that  $f\le z(\delta)$  on $\mathcal{I}.$  
Let $\Omega(\mathcal{I})$  be the component of $\Omega\setminus\mathcal{I}$  whose closure contains $\mathcal{O}.$  
Then (\ref{wheat}) implies that 
\[
f \le \sup_{\mathcal{I}}f+q(p(\delta)) \le z(\delta)+q(p(\delta)) = z_{2}-2p(\delta) \ \ \ {\rm in} \ \Omega(\mathcal{I})
\] 
and so $\limsup_{E\ni {\bf w}\to {\bf o}} c({\bf w})  \le z_{2}-2p(\delta)<z_{2},$  which is a contradiction; hence no such curve $\mathcal{I}$  exists and 
$\mathcal{O}\in \overline{V(\delta)}.$    

Now $f\ge z(\delta)$  in $V(\delta)$   for each $\delta\in (0,p^{-1}(\tau)).$
Let $\mathcal{C}$  be any curve in $\Omega$  which starts at a point ${\bf x}_{0}\in C'_{\rho(p^{-1}(\tau),{\bf o})}({\bf o})$  and   
ends at $\mathcal{O}$  such that 
\[
\mathcal{C}\subset V(\delta) \ \ \ {\rm for\ each} \ \ \delta\in (0,p^{-1}(\tau)).
\]
Since $\liminf_{\mathcal{C}\ni {\bf x}\to \mathcal{O}} f({\bf x}) \ge \lim_{\delta\downarrow 0} z(\delta)=z_{2}$  and 
$z_{2}= \limsup_{\Omega\ni {\bf x}\to \mathcal{O}} f({\bf x}),$  we see that 
\begin{equation}
\label{Because}
\lim_{\mathcal{C}\ni {\bf x}\to \mathcal{O}} f({\bf x})=z_{2}.
\end{equation}
We may, if we wish, extend $\mathcal{C}$  by adding to $\mathcal{C}$  a curve from ${\bf x}_{0}$  to a point on 
$\partial\Omega\setminus \overline{E'_{\rho(p^{-1}(\tau),{\bf o})}({\bf o})}.$  

Now we modify the argument in the proof of \cite[Theorem 2]{NoraKirk1} to show that $Rf(\theta)=z_{2}$  for all $\theta\in (\alpha,\beta);$ 
that is, we shall show that the nontangential limit of $f$  at $\mathcal{O}$  exists and equals $z_{2}.$  
Let $\alpha',\beta'\in (\alpha,\beta)$  with $\alpha'<\beta'.$   

\begin{figure}[htb]
\centerline{
\includegraphics[width=4in]{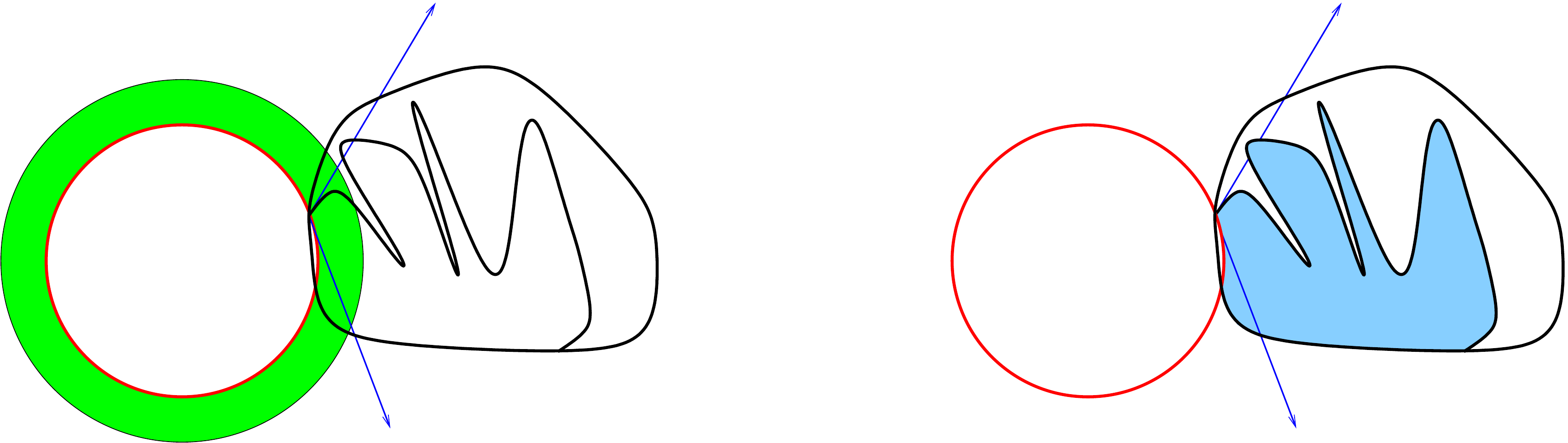}}
\caption{$\Omega,$  $\mathcal{A}_{-}$  and $\mathcal{C}$  (left) \hspace{1cm} $\Omega_{2}$  (right)  \label{Fig1}}
\end{figure}

Let $H_{0}=\sup_{B_{\delta_{0}}(\mathcal{O})\cap\Omega}H$  and fix $c_{0}\in \left(-\frac{1}{4c_{0}H_{0}},0\right).$ 
Set $r_{1}=\frac{1-\sqrt{1+4c_{0}H_{0}}}{2H_{0}}$  and $r_{2}=\frac{1+\sqrt{1+4c_{0}H_{0}}}{2H_{0}}$  (see \cite[p. 171]{LS1}, \cite{FinnMoon}).
Let $\mathcal{A}_{\pm}$  be annuli with inner boundaries $\partial_{1}\mathcal{A}_{\pm}$  
with equal radii $r_{1}$  and outer boundaries $\partial_{2}\mathcal{A}_{\pm}$  with equal radii $r_{2}$  such that 
$\mathcal{O}\in \partial_{1}\mathcal{A}_{+} \cap \partial_{1}\mathcal{A}_{-},$
$\partial_{1}\mathcal{A}_{+}$  is tangent to the ray $\theta=\beta'$  at $\mathcal{O},$  
$\partial_{1}\mathcal{A}_{-}$  is tangent to the ray $\theta=\alpha'$  at $\mathcal{O}$  
and $\partial_{1}\mathcal{A}_{\pm}\cap \{(r\cos\theta,r\sin\theta) : 0<r<\delta_{0}, \alpha'<\theta<\beta'\}=\emptyset$  (see Figure \ref{Fig1}).   
Let $h_{\pm}=h(\hat r_{\pm})$  denote unduloid surfaces defined respectively on $\mathcal{A}_{\pm}$  with constant mean curvature $-H_{0}$  which become 
vertical at $\hat r_{\pm}=r_{1},r_{2}$  and make contact angles of $\pi$  and $0$  with the vertical cylinders $\hat r_{\pm}=r_{2}$  and $\hat r_{\pm}=r_{1}$  
respectively, where $\hat r_{+}({\bf x})=|{\bf x}-{\bf c}_{+}|,$   $\hat r_{-}({\bf x})=|{\bf x}-{\bf c}_{-}|,$  
${\bf c}_{+}$  denotes the center of the annulus $\mathcal{A}_{+}$  and ${\bf c}_{-}$  denotes the center of the annulus $\mathcal{A}_{-}.$  
With respect to the upward direction, the graphs of $h_{\pm}$  over $\mathcal{A}_{\pm}$   have constant mean curvature $-H_{0}$  and 
the graphs of $-h_{\pm}$  over $\mathcal{A}_{\pm}$   have constant mean curvature $H_{0}.$  

\begin{figure}[htb]
\centerline{
\includegraphics[width=2.5in]{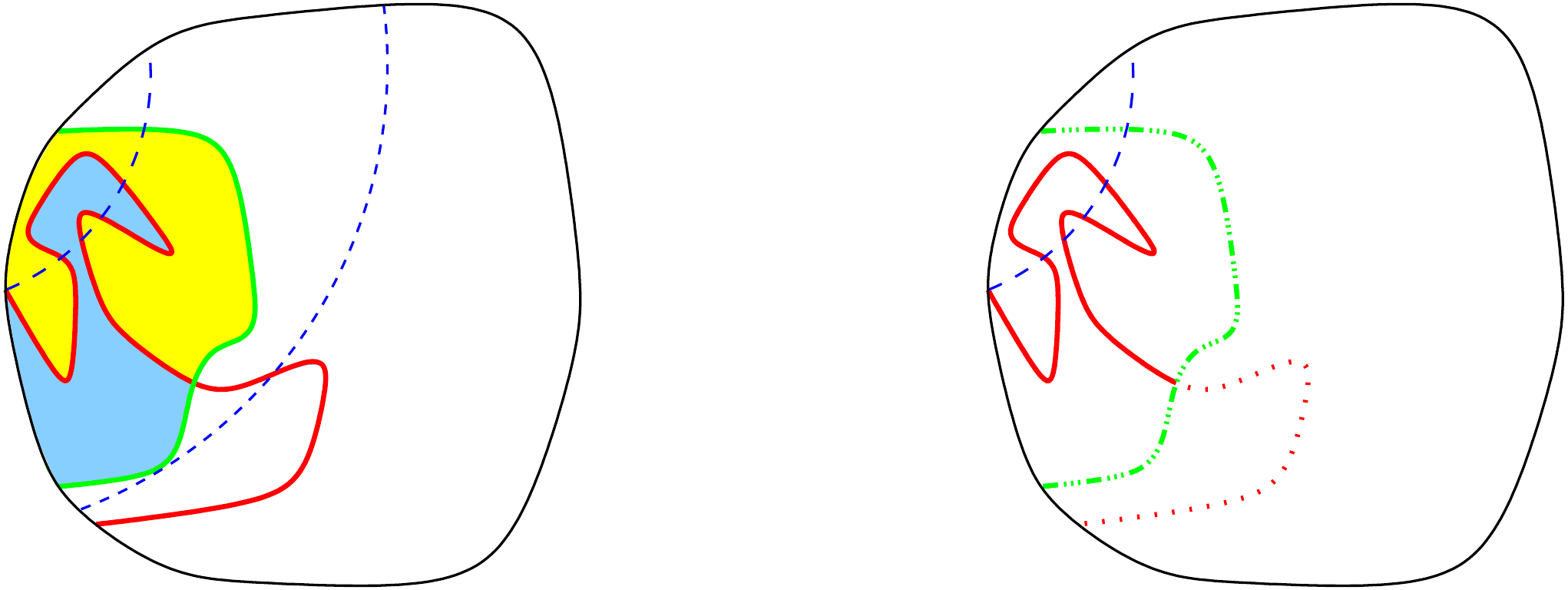}}
\caption{Left: $U_{+}$  (yellow),  $U_{-}$  (blue); Right: $\mathcal{C}(\delta)$ (red) \label{CCC}}
\end{figure}

Set $\tau_{1}=\min\{\tau,r_{2}-r_{1}\}.$  
Let $\delta\in (0,p^{-1}(\tau_{1})).$  Since $\mathcal{C}$  is a curve in $\Omega$  with $\mathcal{O}$  as an endpont, there exists 
${\bf x}(\delta)\in \mathcal{C}\cap C_{\rho(\delta,{\bf o})}'({\bf o})$  such that the portion $\mathcal{C}(\delta)$  of $\mathcal{C}$  between 
$\mathcal{O}$  and ${\bf x}(\delta)$  lies in $E'_{\rho(\delta,{\bf o})}({\bf o})$  and divides $E'_{\rho(\delta,{\bf o})}({\bf o})$  into two components. 
Let $U_{+}$  be the component of $E'_{\rho(\delta,{\bf o})}({\bf o}) \setminus \mathcal{C}(\delta)$  whose closure contains a portion of $\partial^{+}\Omega$
and  $U_{-}$  be the component of $E'_{\rho(\delta,{\bf o})}({\bf o}) \setminus \mathcal{C}(\delta)$  whose closure contains a portion of $\partial^{-}\Omega$ 
(see Figure \ref{CCC} with $C_{\rho(\delta,{\bf o})}'({\bf o})$  (green) and $\mathcal{C}$  (red)).

\begin{figure}[htb]
\centerline{
\includegraphics[width=2in]{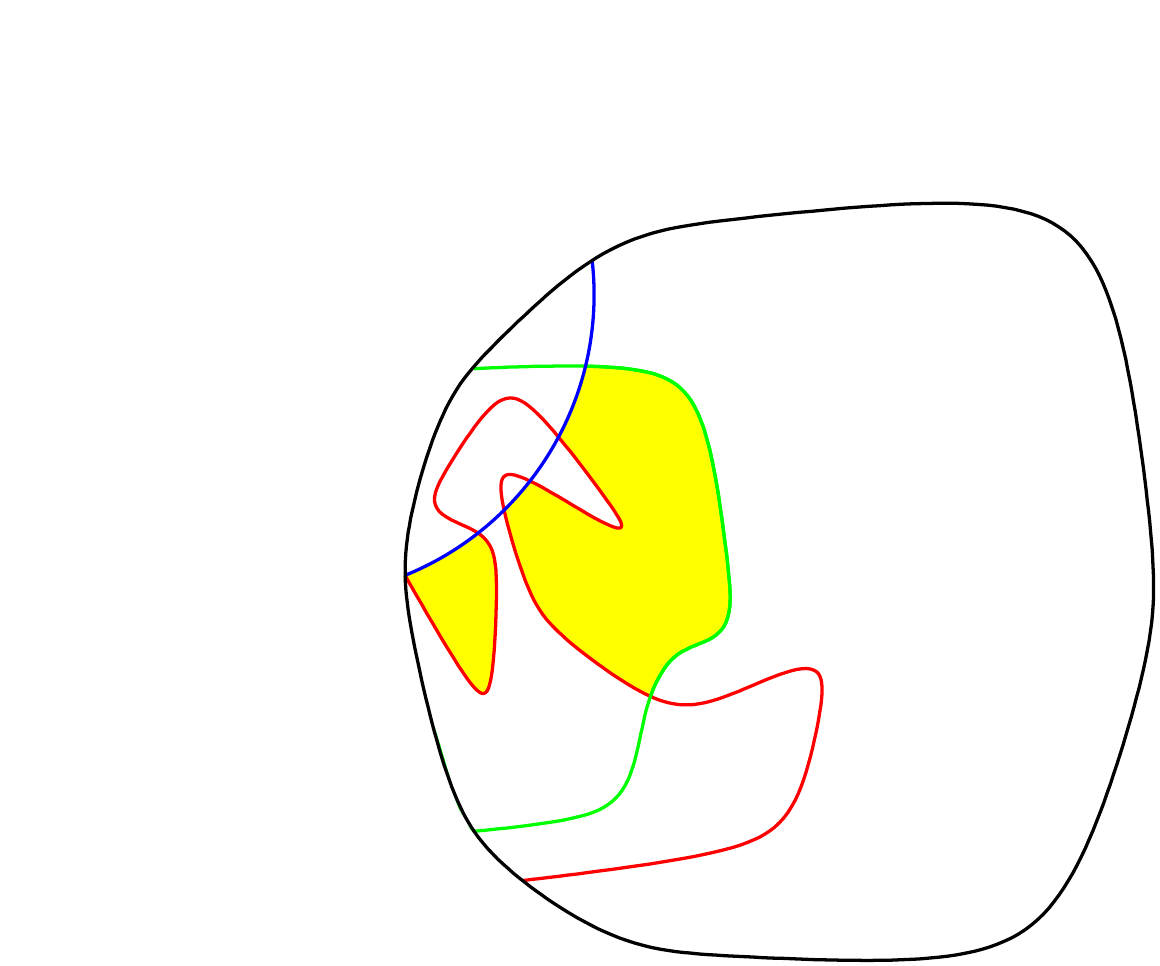}}
\caption{ $\Omega$  and  $\Omega_{1}$   \label{Fig2}}
\end{figure}

Since $\mathcal{C}(\delta)\subset V(\delta),$
\[
f({\bf x})\ge z(\delta) \ \ \ \ {\rm for}\ \ {\bf x}\in \mathcal{C}(\delta)
\]
and, in particular, $f({\bf x}(\delta))\ge z(\delta).$
Since $|f({\bf x}(\delta))-f({\bf y})|\le l_{\rho(\delta,{\bf o})}<p(\delta)$  for ${\bf y}\in C'_{\rho(\delta,{\bf o})}({\bf o}),$  we see that 
\[
f \ge z(\delta)-p(\delta) \ \ \ \ {\rm on} \ \  C'_{\rho(\delta,{\bf o})}({\bf o}) \cup \mathcal{C}(\delta).
\]
Let $q_{2}$  denote a modulus of continuity of $-h(\hat r_{+}).$ 
Then 
\[
f \ge z(\delta)-p(\delta)-q_{2}(p(\delta)) \ \ \ \ {\rm in} \ \ U_{+}\setminus \overline{B_{r_{1}}(c_{+})}. 
\]
Thus 
\[
\liminf_{U_{+}\setminus \overline{B_{r_{1}}(c_{+})} \ni {\bf x}\to\mathcal{O}}f({\bf x})\ge z_{2}.
\]
If we set $\Omega_{1}=U_{+}\setminus \overline{B_{r_{1}}({\bf c}_{+})}$  and recall that $z_{2}= \limsup_{\Omega\ni {\bf x}\to \mathcal{O}} f({\bf x}),$  
we have  
\begin{equation}
\label{AAA}
\lim_{\Omega_{1}\ni {\bf x}\to \mathcal{O}} f({\bf x})=z_{2}.
\end{equation}
(We note that $\Omega_{1}$  might not be connected (see Figure \ref{Fig2}) and might even have an infinite number of components but one sees that this does not affect 
the comparison argument which establishes (\ref{AAA}).)  
In a similar manner, we see that 
\begin{equation}
\label{BBB}
\lim_{\Omega_{2}\ni {\bf x}\to \mathcal{O}} f({\bf x})=z_{2},
\end{equation}
where $\Omega_{2}=U_{-}\setminus \overline{B_{r_{1}}({\bf c}_{-})}.$ 
Since $\Omega_{1}\cup\Omega_{2}\cup \mathcal{C}(\delta)=E'_{\rho(\delta,{\bf o})}({\bf o})
\setminus\left(\overline{B_{r_{1}}({\bf c}_{+})\cup B_{r_{1}}({\bf c}_{-})}\right),$
we see that $Rf(\theta)=z_{2}$  for each $\theta\in (\alpha',\beta').$  
Since $\alpha'$  and $\beta'$  are arbitrary (with $\alpha<\alpha'<\beta'<\beta$), Theorem \ref{Thm1} is proven.  \qed

 \begin{rem}
 \label{Why}
 If $\beta-\alpha<\pi,$  then the existence of the one-sided barrier $\psi^{+}$  is uncertain and so the proof that a curve (i.e. $\mathcal{C}$) in $\Omega$ 
 with $\mathcal{O}$  as an endpoint such that (\ref{Because}) holds is uncertain.  
 \end{rem}
\vspace{3mm}

\noindent {\bf Proof of Theorem \ref{Thm2}:}  
All of the claims in the theorem except those in the last sentence follow from \cite[Theorem 1]{NoraKirk1} and \cite{EchartLancaster1} (when $\beta-\alpha > \pi$) 
and \cite[Theorem 2]{NoraKirk1}  and \cite{EchartLancaster1} (when $\beta-\alpha < \pi$).  
(When $\beta-\alpha = \pi,$  all of the claims follow from Theorem \ref{Thm1}  and \cite{NoraKirk3}.) 
The claims follow once we prove that the results of \cite{NoraKirk3} hold under the assumptions of Theorem \ref{Thm2}. 
Let us assume 
\begin{equation}
\label{curve3}
\limsup_{\partial^{\pm}\Omega\ni {\bf x}\to \mathcal{O}}\left(\kappa({\bf x})+2|H({\bf x})|\right)<0.
\end{equation}

Suppose $\beta-\alpha > \pi.$   Let $\delta_{1}>0$  be small enough that $B_{\delta_{1}}(\mathcal{O})\cap \Omega\setminus \Sigma^{+}$  has two components.  
Let $\Omega_{+}$  be the component whose closure contains $B_{\delta_{1}}(\mathcal{O})\cap \partial^{+}\Omega$  and notice that the tangent directions to 
$\partial\Omega_{+}$  at $\mathcal{O}$  are $\alpha'=\beta-\pi$  and $\beta$  and the curvature $\kappa_{+}(\mathcal{O})$  of $\partial\Omega_{+}$   at $\mathcal{O}$ 
satisfies 
\[
\kappa_{+}(\mathcal{O})<-2H(\mathcal{O})
\]
since $\kappa_{+}({\bf x})=\kappa({\bf x})$  for ${\bf x}\in B_{\delta_{1}}(\mathcal{O})\cap \partial^{+}\Omega$  and (\ref{curve3})  implies
\[
\kappa_{+}(\mathcal{O})=\limsup_{\partial^{+}\Omega\ni {\bf x}\to \mathcal{O}} \kappa_{+}({\bf x})<-2|H(\mathcal{O})|.
\]
By restricting $f$  to $\Omega_{+},$  we see that the existence of $Rf(\beta)$  follows from \cite{NoraKirk3}. 
A similar argument implies $Rf(\alpha)$  also exists.

Suppose $\beta-\alpha < \pi.$  Then $Rf(\alpha)$  exists and equals $m.$  
Let $\delta_{1}>0$  be small enough that $B_{\delta_{1}}(\mathcal{O})\setminus \Sigma^{+}$  has two components and let $\Omega_{+}$  be the component 
which contains $B_{\delta_{1}}(\mathcal{O})\cap \Omega.$  
Then the tangent directions to $\partial\Omega_{+}$  at $\mathcal{O}$  are $\alpha'=\beta-\pi$  and $\beta$  and, as before, the curvature $\kappa_{+}(\mathcal{O})$  
of $\partial\Omega_{+}$   at $\mathcal{O}$  satisfies $\kappa_{+}(\mathcal{O})<-2H(\mathcal{O}).$  
Thus upper and lower Bernstein pairs $(U^{\pm},\psi^{\pm})$  exist for $\Gamma=B_{\delta_{2}}(\mathcal{O}) \cap \partial\Omega_{+}$  and $H$  when 
$\delta_{2}\in (0,\delta_{1})$  is sufficiently small and $U^{\pm}=B_{\delta_{2}}(\mathcal{O}) \cap\Omega_{+}.$
We may parametrize $S_{1}=S_{0}\cap \left(B_{\delta_{2}}(\mathcal{O})\times\Real\right)$  in isothermal coordinates 
\begin{equation}
\label{PARAMETRIC}
Y(u,v) = (a(u,v),b(u,v),c(u,v)) \in C^{2}(E:S_{1}) 
\end{equation}
as in \cite{NoraKirk3}  with the properties noted there (e.g. $a_{1},\dots,a_{5}$) and prove in essentially the same manner as in \cite{NoraKirk3} that 
$Y$  is uniformly continuous on $E$  and so extends to a continuous function on $\overline{E}.$  
(Notice the similarity of methods used in \cite{NoraKirk1} and \cite{NoraKirk3}.) 
The existence of $Rf(\beta)$  then follows as in \cite{NoraKirk3}.  \qed

\begin{figure}[htb]
\centerline{
\includegraphics[width=4in]{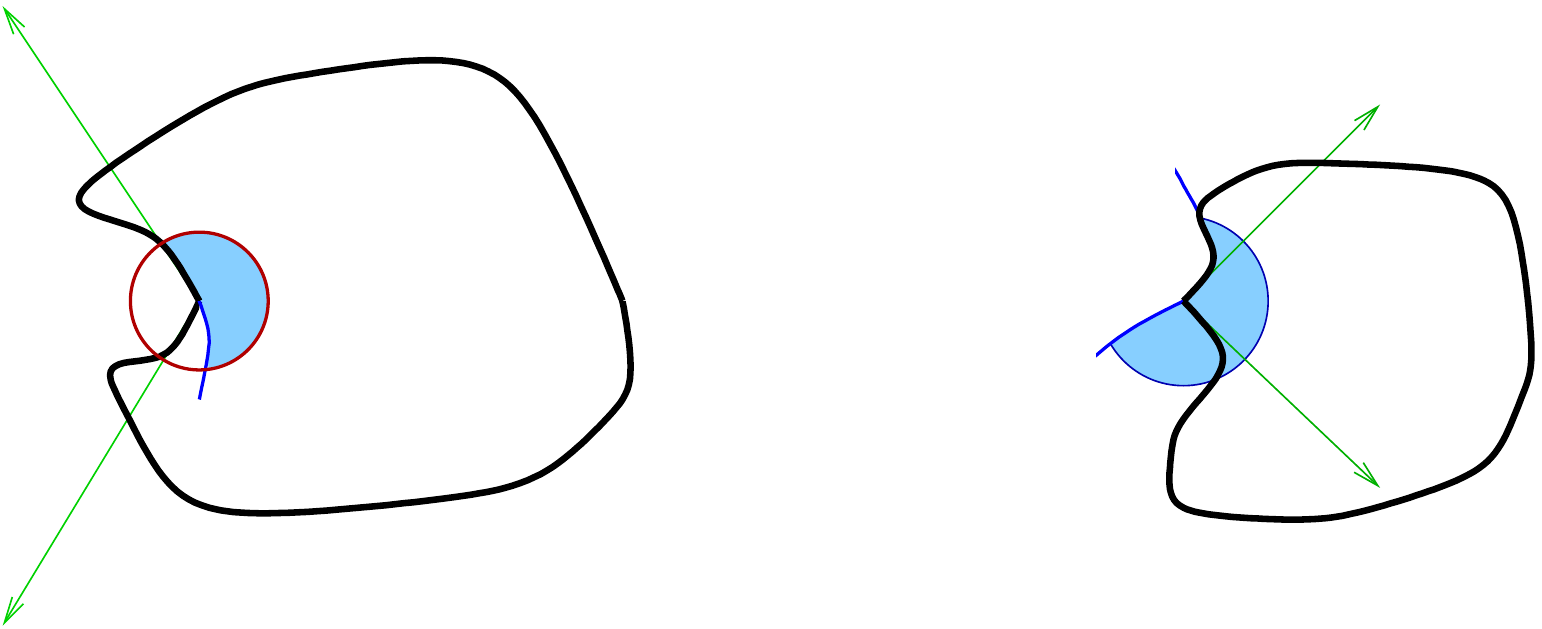}}
\caption{$\Omega_{+}$  when $\beta-\alpha > \pi$  (left) \hspace{3mm} $\Omega_{+}$  when $\beta-\alpha < \pi$  (right)  \label{Fig3}}
\end{figure}

\end{document}